\documentclass[11pt]{article}
\usepackage{amssymb,amsmath,latexsym}
\usepackage{tablists}
\usepackage{hyperref}
\usepackage[latin1]{inputenc}
\usepackage[french,english]{babel}
\usepackage{graphicx}
\newtheorem{ex}{Example}
\newtheorem{qu}{Question}
\newtheorem{remark}{\\Remark}

\newtheorem{lemma}{lemma}
\newtheorem{thh}{Theorem }
\newtheorem{co}{Corollary}
\newcommand{\p}{{{\bf Proof.\,\,}}}

\newcommand{\biindice}[3]%
{

\begin{array}[t]{c}
#1\\
{\scriptstyle #2}\\
{\scriptstyle #3}
\end{array}

}

\begin{document}
\title{\Large  \textbf{Problem of descent spectrum equality }}
\author{\normalsize Abdelaziz Tajmouati$^{1}$ and Hamid Boua$^{2}$\\\\
\normalsize Sidi Mohamed Ben Abdellah University\\
\normalsize Faculty of Sciences Dhar El Mahraz \\
\normalsize Laboratory of Mathematical Analysis and Applications\\
\normalsize Fez, Morocco\\
\normalsize Email: abdelaziz.tajmouati@usmba.ac.ma\,\,\,\, hamid12boua@yahoo.com
}

\date{}\maketitle
\begin{abstract}
Let $\mathcal{B}(X)$ be the algebra of all bounded operators acting on an infinite dimensional complex Banach space $X$.
We say that an operator $T \in \mathcal{B}(X)$ satisfies the \textbf{problem of descent spectrum equality }, if the descent spectrum of $T$ as an operator coincides with the descent spectrum of $T$ as an element of the algebra of all bounded linear operators on $X$. In this paper we are interested in the problem of descent spectrum equality . Specifically, the problem is to consider the following question: Let  $T \in \mathcal{B}(X)$ such that $\sigma(T)$ has non empty interior,  under which condition on $T$ does $\sigma_{desc}(T)=\sigma_{desc}(T, \mathcal{B}(X))$ ?
\end{abstract}
\section{\large Introduction}
In this paper, $X$ denotes a complex Banach space and $\mathcal{B}(X)$ denotes the Banach algebra of all bounded linear
operators on $X$. Let $T\in \mathcal{B}(X)$, we denote by $R(T)$, $N(T)$, $\rho(T)$, $\sigma(T)$, $\sigma_{p}(T)$, $\sigma_{ap}(T)$ and $\sigma_{su}(T)$ respectively the range, the kernel, the resolvent set, the spectrum, the point spectrum, the approximate point spectrum and the surjectivity spectrum of $T$. It is well known that $\sigma(T)=\sigma_{su}(T)\cup \sigma_{p}(T)=\sigma_{su}(T)\cup \sigma_{ap}(T)$.
The ascent of $T$ is defined by $a(T)=\min\{p : N(T^p)=N(T^{p+1})\}$, if no such $p$ exists, we let $a(T)=\infty$.
Similarly, the descent of $T$ is $d(T)=\min\{q :  R(T^q)=R(T^{q+1})\}$, if no such $q$ exists, we let $d(T)=\infty$ \cite{Aie}, \cite{lay} and \cite{Mul}. It is well known that if both $a(T)$ and $d(T)$ are finite then $a(T)=d(T)$ and we have the decomposition $X=R(T^p)\oplus N(T^p)$ where $p=a(T)=d(T)$. The  descend  and ascent  spectrum are defined by:
$$\sigma_{desc}(T)=\{\lambda \in \mathbb{C}  :  d(\lambda -T)=\infty \}$$
$$\sigma_{asc}(T)=\{\lambda \in \mathbb{C}  :  a(\lambda -T)=\infty \}$$
$\mathcal{A}$ will denote a complex Banach algebra with unit. For every $a \in \mathcal{A}$, the left multiplication operator $L_a$ is given by $L_a(x) = ax$ for all $x\in \mathcal{A}$. By definition the descent of an element $a\in \mathcal{A}$ is $d(a) := d(L_a)$, and the descent spectrum of $a$ is the set
$\sigma_{desc}(a) := \{\lambda \in \mathbb{C} : d(a-\lambda) = \infty \}.$\\In general $\sigma_{desc}(T)\subseteq \sigma_{desc}(T, \mathcal{B}(X))$, and we say that an operator $T$ satisfies the descent spectrum equality whenever, the descent spectrum of $T$ as an operator coincides with the descent spectrum of $T$ as an element of the algebra of all bounded linear operators on $X$.\\
The operator $T\in \mathcal{B}(X)$ is said to have the single-valued extension property at $\lambda_0 \in \mathbb{C}$,
abbreviated $T$ has the SVEP at $\lambda_0$, if for every neighbourhood $\mathcal{U}$ of $\lambda_0$ the
only analytic function $f : \mathcal{U} \rightarrow X$ which satisfies the equation $(\lambda I - T)f(\lambda) = 0$
is the constant function $f \equiv 0$. For an arbitrary operator $T \in B(X)$ let $\mathcal{S}(T)= \{\lambda \in \mathbb{C}   :  T \mbox{ does not have the SVEP at } \lambda \}$. Note that $\mathcal{S}(T)$ is open and is contained in the interior of the point spectrum $\sigma_p(T)$.
The operator $T$ is said to have the SVEP if $\mathcal{S}(T)$ is empty. According to \cite{Lar} we have $\sigma(T) = \sigma_{su}(T)\cup  \mathcal{S}(T)$.\\
For an operator $T \in \mathcal{B}(X)$ we shall denote by $\alpha(T)$ the dimension of the kernel $N(T)$, and by $\beta(T)$ the codimension of the range $R(T)$.  We recall that an operator $T \in \mathcal{B}(X)$ is called upper semi-Fredholm if $\alpha(T)< \infty$  and $R(T)$ is closed,
while $T \in \mathcal{B}(X)$ is called lower semi-Fredholm if $\beta(T)< \infty$. Let $\Phi_{+}(X)$ and $\Phi_{-}(X)$ denote the
class of all upper semi-Fredholm operators and the class of all lower semi-Fredholm operators,
respectively. The class of all semi-Fredholm operators is defined by $\Phi_{±}(X) := \Phi_{+}(X)\cup \Phi_{-}(X)$, while the class of all Fredholm operators is defined by $\Phi(X) := \Phi_{+}(X)\cap \Phi_{-}(X)$. If $T \in \Phi_{±}(X)$, the index of $T$ is defined by ind$(T) := \alpha(T) - \beta(T)$. The class of all upper semi-Browder operators is defined by $B_{+}(X) := \{T \in \Phi_{+}(X) : a(T) < \infty \}$, the upper semi-Browder spectrum of $T \in \mathcal{B}(X)$ is defined by $\sigma_{ub}(T) := \{\lambda \in \mathbb{C}: \lambda I - T \notin B_{+}(X)\}$. The class of all upper semi-Weyl operators is defined by $W_{+}(X) := \{T \in \Phi_{+}(X) : \mbox{ind}(T) \leq 0 \}$, the upper semi-Weyl spectrum is defined by $\sigma_{uw}(T) := \{\lambda \in \mathbb{C}: \lambda I - T \notin W_{+}(X)\}$. \\
Recently, Haily, Kaidi and Rodrigues Palacios \cite{H} have studied and characterized the Banach spaces verifying property descent spectrum equality, (Banach which are isomorphic to $\ell^1(I)$ or $\ell^2(I)$ for some set $I$, the not isomorphic to any of its proper quotients...). On the other hand, they have shown that if $T \in \mathcal{B}(X)$ with a spectrum $\sigma(T)$ of empty interior, then $\sigma_{desc}(T)=\sigma_{desc}(T, \mathcal{B}(X))$.\\It is easy to construct an operator $T$ satisfying the descent spectrum equality such that the interior of the point spectrum $\sigma(T)$ is nonempty. For example, let $T$ the bilateral right shift on the Hilbert space $\ell^2(\mathbb{Z})$, so that $T(x_{n})_n=(x_{n-1})_n$ for all $(x_{n})_{n\in \mathbb{Z}}\in \ell^2(\mathbb{Z})$. It is easily seen that $\sigma(T)=\overline{\mathbb{D}}$ closed unit disk. Since $\ell^2(\mathbb{Z})$ is a Hilbert space, then $\sigma_{desc}(T)=\sigma_{desc}(T, \mathcal{B}(X))$. Motivated by the previous Example, our goal is to study the following question:
\begin{qu}
Let $T \in \mathcal{B}(X)$. If $\sigma(T)$ has non empty interior,  under which condition on $T$ does
$\sigma_{desc}(T)=\sigma_{desc}(T, \mathcal{B}(X))$ ?
\end{qu}

\section{Main results}
We start by the following lemmas.
\begin{lemma}\cite{H}
Let $T$ be in $\mathcal{B}(X)$ with finite descent $d=d(T)$. Then there exists $\delta > 0$ such that, for every $\mu \in \mathbb{K}$ with $0 < |\mu| < \delta$, we have:
\begin{enumerate}
\item $T-\mu$ is surjective,
\item $\mbox{dim} N(T-\mu) =\mbox{dim}(N(T) \cap R(T^d))$.
\end{enumerate}
\end{lemma}

\begin{lemma}\cite{H}
 Let $X$, $Y$ and $Z$ be Banach spaces, and let $F : X \rightarrow Z$ and $G : Y \rightarrow Z$ be bounded linear operators such that $N(G)$ is complemented in $Y$, and $R(F)\subseteq R(G)$. Then there exists a bounded linear operator $S : X \rightarrow Y$
satisfying $F = GS$.
\end{lemma}

We have the following theorem.
\begin{thh}\label{a}
Let $T \in \mathcal{B}(X)$ and $D\subseteq \mathbb{C}$ be a closed subset such that $\sigma(T)=\sigma_{su}(T)\cup D$, then
\begin{center}
$\sigma_{desc}(T)\cup \mbox{int}(D) = \sigma_{desc}(T, \mathcal{B}(X))\cup \mbox{int}(D)$
\end{center}
\end{thh}
\p
Let $\lambda$ be a complex number such that $T-\lambda$ has finite descent $d$ and $\lambda \notin int(D)$. According to lemma 1, there is $\delta > 0$ such that, for every $\mu \in \mathbb{C}$ with $0 < |\lambda - \mu| < \delta$, the operator $T - \mu$ is surjective and $\dim N(T-\mu)=\dim N(T-\lambda)\cap R(T-\lambda)^d$. Let $D^{*}(\lambda, \delta)=\{\mu\in \mathbb{C} : 0 < |\lambda - \mu| < \delta\}$. Since $\lambda \notin int(D)$, then $D(\lambda, \delta)\backslash D \neq \emptyset $ is  non-empty open subset of $\mathbb{C}$. Let $\lambda_0 \in D^{*}(\lambda, \delta)\backslash D $, then $T-\lambda_0$ is invertible, hence the continuity of the index ensures that ind$(T - \mu) = 0$ for all $\mu \in D^{*}(\lambda, \delta)$. But for $\mu \in D^{*}(\lambda, \delta), T - \mu $ is surjective, so it follows that $T - \mu$ is invertible. Therefore, $\lambda$ is isolated in $\sigma(T)$. By \cite[Theorem 3.81]{Aie}, we have $\lambda$ is a pole of the resolvent of $T$. Using \cite[Theorem V.10.1]{lay}, we obtain  $T - \lambda$ has a finite descent and a finite ascent and $X = N((T - \lambda)^d) \oplus R((T - \lambda)^d)$. It follows that $N((T - \lambda)^d)$ is complemented in $X$. Applying lemma 2, there exists $S \in \mathcal{B}(X)$ satisfying $(T -\lambda)^d = (T -\lambda)^{d+1}S$, which forces that $\lambda \notin \sigma_{desc}(T, \mathcal{B}(X))\cup int(D)$.

\begin{co}
Let $T\in \mathcal{B}(X)$. If $T$ satisfies any of the conditions following:
\begin{enumerate}
\item $\sigma(T)=\sigma_{su}(T)$,
\item int($\sigma_{ap}(T))=\emptyset$,
\item int($\sigma_{p}(T))=\emptyset$,
\item int($\sigma_{asc}(T))=\emptyset$,
\item int($\sigma_{ub}(T))=\emptyset$,
\item int($\sigma_{uw}(T))=\emptyset$,
\item $ \mathcal{S}(T)=\emptyset$.
\end{enumerate}
Then
\begin{center}
$\sigma_{desc}(T)=\sigma_{desc}(T, \mathcal{B}(X)) $
\end{center}
\end{co}
\p

The assertions 1, 2, 3, and 7 are obvious.\\
4. Note that, $\sigma(T)=\sigma_{su}(T)\cup \sigma_{asc}(T)$. Indeed, let $\lambda \notin \sigma_{su}(T)\cup \sigma_{asc}(T)$, then $T-\lambda$ is surjective and $T-\lambda$ has finite ascent, therefore $a(T-\lambda)=d(T-\lambda)= 0$, and hence $\lambda \notin \sigma(T)$. If int($\sigma_{asc}(T))=\emptyset$, by theorem 1, we have $\sigma_{desc}(T)=\sigma_{desc}(T, \mathcal{B}(X)) $.\\
5. If int($\sigma_{ub}(T))=\emptyset$, then int($\sigma_{asc}(T))=\emptyset$,  therefore $\sigma_{desc}(T)=\sigma_{desc}(T, \mathcal{B}(X)) $.\\
6. Note that, $\sigma(T)=\sigma_{su}(T)\cup \sigma_{uw}(T)$. Indeed, let $\lambda \notin \sigma_{su}(T)\cup \sigma_{uw}(T)$, then $T-\lambda$ is surjective and ind$(T-\lambda)\leq 0$, therefore ind$(T-\lambda)=\dim N(T-\lambda)= 0$, and
hence $\lambda \notin \sigma(T)$. If int($\sigma_{uw}(T))=\emptyset$, by theorem 1, we have $\sigma_{desc}(T)=\sigma_{desc}(T, \mathcal{B}(X)) $.

\begin{ex}
We consider the Césaro operator $C_p$ defined on the classical Hardy space $\mathbb{H}_p(\mathbb{D})$, $\mathbb{D}$ the open unit disc
and $1 < p < \infty$. The operator $C_p$ is defined by $(C_pf)(\lambda) := \frac{1}{\lambda}\int_0^{\lambda} \frac{f(\mu)}{1-\mu}d \mu$ or all $f \in \mathbb{H}_p(\mathbb{D})$ and $\lambda \in \mathbb{D}$.
As noted by T.L. Miller, V.G. Miller and Smith \rm\cite{TM}, the spectrum
of the operator $C_p$ is the entire closed disc $\Gamma_p$, centered at $p/2$ with radius
$p/2$, and $\sigma_{ap}(C_p)$ is the boundary $\partial \Gamma_p$, then  int($\sigma_{ap}(C_p)$)=int($\sigma_{p}(C_p))=\emptyset$. By applying corollary 1, then $\sigma_{desc}(C_p)=\sigma_{desc}(C_p,\mathcal{B}(\mathbb{H}_p(\mathbb{D}))) $.
\end{ex}

\begin{ex}
Suppose that $T$ is an unilateral weighted right shift on $\ell^p(\mathbb{N})$, $1 \leq p < \infty$, with weight sequence  $(\omega_n)_{n\in \mathbb{N}}$, $T$ is the operator defined by: $Tx:=\sum_{n=1}^{\infty}\omega_nx_ne_{n+1}$ for all $x:=(x_n)_{n\in \mathbb{N}}\in \ell^p(\mathbb{N})$
. If
$c(T) = \lim_{n \rightarrow +\infty} \inf(\omega_1...\omega_n)^{1/n} = 0$, by \rm\cite[Corollary 3.118]{Aie}, we have $T$ has SVEP. By applying corollary 1, then $\sigma_{desc}(T)=\sigma_{desc}(T,\mathcal{B}(X)) $.
\end{ex}

A mapping $T : \mathcal{A} \rightarrow \mathcal{A}$ on a commutative complex Banach algebra $\mathcal{A}$ is said to
be a multiplier if:
\begin{center}
$u(Tv) = (Tu)v \mbox{ for all } u, v \in \mathcal{A}.$
\end{center}
Any element $a \in \mathcal{A}$ provides an example, since, if $L_a : \mathcal{A} \rightarrow \mathcal{A}$ denotes the
mapping given by $L_a(u) := au$ for all $u \in \mathcal{A}$, then the multiplication operator La
is clearly a multiplier on $\mathcal{A}$. The set of all multipliers of $\mathcal{A}$ is denoted by $M(\mathcal{A})$. We recall that an algebra $\mathcal{A}$ is said to be semi-prime if $\{0\}$ is the only
two-sided ideal $J$ for which $J^2 = {0}$.
\begin{co}
Let $T \in M(\mathcal{A})$ be a multiplier on a semi-prime commutative Banach algebra $\mathcal{A}$
then:
\begin{center}
$\sigma_{desc}(T)=\sigma_{desc}(T,\mathcal{B}(X)) $
\end{center}
\end{co}

\p
If $T \in M(\mathcal{A})$, from \cite[Proposition 4.2.1]{Aie}, we have $\sigma(T)=\sigma_{su}(T)$. By applying corollary 1, then: $\sigma_{desc}(T)=\sigma_{desc}(T,\mathcal{B}(X)) $.

\begin{thh}
Let $T \in \mathcal{B}(X)$. If for every connected component $G$ of $\rho_{desc}(T)$ we have that $G\cap \rho(T)\neq \emptyset $, then
\begin{center}
$\sigma_{desc}(T)=\sigma_{desc}(T, \mathcal{B}(X))$
\end{center}
\end{thh}

\p
Let $\lambda$ be a complex number such that $T-\lambda$ has finite descent $d$. According to lemma 1, there is $\delta > 0$ such that, for every $\mu \in \mathbb{C}$ with $0 < |\lambda - \mu| < \delta$, the operator $T - \mu$ is surjective and $\dim N(T-\mu)=\dim N(T-\lambda)\cap R(T-\lambda)^d$. $D^{*}(\lambda, \delta)=\{\mu\in \mathbb{C} : 0 < |\lambda - \mu| < \delta\}$ is a connected subset of $\rho_{desc}(T)$, then there exists a connected component  $G$ of $\rho_{desc}(T)$ contains $D^{*}(\lambda, \delta)$. Since $ G\cap \rho(T)$ is non-empty hence the continuity of the index ensures that $\mbox{ind}(T-\mu)=0$ for all $\mu \in D^{*}(\lambda, \delta)$. But for $\mu \in G$, $T - \mu$ is surjective, so it follows that $T - \mu$ is invertible.
Thus $G\subseteq \rho(T)$, therefore, $\lambda$ is isolated in $\sigma(T)$. Consequently $\lambda \notin \sigma_{desc}(T, \mathcal{B}(X))$, which completes the proof.

\begin{remark}
We recall that an operator $R \in \mathcal{B}(X)$ is said to be Riesz if $R-\lambda$ is Fredholm for every non-zero complex number $\lambda$. From \rm\cite{lay}, $\sigma_{desc}(R)=\{0\}$, then for every connected component $G$ of $\rho_{desc}(R)$, we have that $G\cap \rho(R)\neq \emptyset $. Consequently $\sigma_{desc}(R, \mathcal{B}(X))=\{0\}$
\end{remark}

\begin{ex}
Consider the unilateral right shift operator $T$ on the space $X :=\ell^p$ for some $1 \leq p \leq \infty$. Because $\sigma(T)=\sigma_{desc}(T)$, then for every $G$ is a connected component of $\rho_{desc}(T)$ we have that $G\cap \rho(T)\neq \emptyset $. Consequently $\sigma_{desc}(T, \mathcal{B}(X))=\overline{\mathbb{D}}$ closed unit disk.
\end{ex}

\begin{thh}
Let $T \in \mathcal{B}(X)$. If for every connected component $G$ of $\rho_{su}(T)$ we have that $G\cap \rho_p(T)\neq \emptyset $, then:
\begin{center}
$\sigma_{desc}(T)=\sigma_{desc}(T, \mathcal{B}(X))$
\end{center}
\end{thh}
\p
Let $\lambda$ be a complex number such that $T-\lambda$ has finite descent $d$. According to lemma 1, there is $\delta > 0$ such that, for every $\mu \in \mathbb{C}$ with $0 < |\lambda - \mu| < \delta$, the operator $T - \mu$ is surjective and $\dim N(T-\mu)=\dim N(T-\lambda)\cap R(T-\lambda)^d$. Therefore $D^{*}(\lambda, \delta)=\{\mu\in \mathbb{C} : 0 < |\lambda - \mu| < \delta\}$ is a connected subset of $\rho_{su}(T)$, then there exists a connected component  $G$ of $\rho_{su}(T)$ contains $D^{*}(\lambda, \delta)$. Since $ G\cap \rho_p(T)$ is non-empty hence the continuity of the index ensures that $\mbox{ind}(T-\mu)=0$ for all $\mu \in D^{*}(\lambda, \delta)$. But for $\mu \in G$, $T - \mu$ is surjective, so it follows that $T - \mu$ is invertible, therefore, $\lambda$ is isolated in $\sigma(T)$. Consequently $\lambda \notin \sigma_{desc}(T, \mathcal{B}(X))$.

\begin{remark}
Let $T \in \mathcal{B}(X)$ an operator such that $\sigma(T)=\sigma_{su}(T)$, then for every connected component $G$ of $\rho_{su}(T)$, we have $G\cap \rho_p(T)\neq \emptyset $. Using Theorem 3, we obtain $\sigma_{desc}(T)=\sigma_{desc}(T, \mathcal{B}(X))$.
\end{remark}

\end{document}